\documentclass[twoside,12pt]{article}
\usepackage[]{amsmath,amsthm,amssymb}
\usepackage[latin1]{inputenc}

\begin{document}

\title{{\Large\bf  Discrete  index  transformations with  Bessel and Lommel functions}}

\author{Semyon  YAKUBOVICH}
\maketitle

\markboth{\rm \centerline{ Semyon   YAKUBOVICH}}{}
\markright{\rm \centerline{  DISCRETE TRANSFORMS WITH BESSEL AND LOMMEL FUNCTIONS}}

\begin{abstract} {\noindent Discrete analogs of the index  transforms, involving Bessel and Lommel functions are introduced and investigated. The corresponding inversion theorems for   suitable classes  of functions and sequences   are established. }

\end{abstract}
\vspace{4mm}



\vspace{4mm}

\section {Introduction and preliminary results}

Our goal in this paper is to investigate the mapping properties and prove inversion formulas for the following transformations between suitable sequences $\{a_n\}_{n\ge 1}$ and functions $f$ in terms of the series and integrals,  which are associated with  Bessel and Lommel functions  (cf. [2],  Ch. 10, 11), namely,

$$f(x)=  \sum_{n=0}^\infty {a_n \over \cosh(\pi n/2)} \ {\rm Re} \left[ J_{in}(x) \right],\quad x > 0,\eqno(1.1)$$

$$ a_n=  {1\over \cosh(\pi n/2)} \int_0^\infty   {\rm Re} \left[ J_{in}(x) \right] f(x) dx,\quad n \in \mathbb{N}_0,\eqno(1.2)$$

$$f(x)=\sum_{n=1}^\infty {a_n  \over \sinh(\pi n/2)} \  {\rm Im} \left[ J_{in}(x) \right],\quad x > 0,\eqno(1.3)$$

$$ a_n=  {1\over \sinh (\pi n/2)} \int_0^\infty   {\rm Im} \left[ J_{in} (x) \right] f(x) dx,\quad n \in \mathbb{N},\eqno(1.4)$$

$$f(x)= \sum_{n=1}^\infty a_n \  \Gamma\left({1-\mu-in\over 2}\right) \Gamma\left({1-\mu+in\over 2}\right) S_{\mu,in} (x),\quad x > 0,\eqno(1.5)$$

$$ a_n=    \Gamma\left({1-\mu-in\over 2}\right) \Gamma\left({1-\mu+in\over 2}\right)   \int_0^\infty  S_{\mu,in} (x) f(x) dx,\quad n \in \mathbb{N}.\eqno(1.6)$$
Here  $i$ is the imaginary unit and ${\rm Re},\ {\rm Im}$ denote the real and imaginary parts of a complex-valued function.   We call transformations (1.1)-(1.6) the discrete index transforms, comparing them with continuum  analogs (cf. [4]). Bessel functions $J_\nu(z),\ Y_\nu(z),\ z,\nu \in \mathbb{C}$ of the first and second kind, respectively, are solutions of the Bessel differential equation 

$$  z^2{d^2u\over dz^2}  + z{du\over dz} + (z^2- \nu^2)u = 0.\eqno(1.7)$$
These functions have the following asymptotic behavior at infinity and near the origin
$$ J_\nu(z) = \sqrt{2\over \pi z} \cos \left( z- {\pi\over 4} (2\nu+1)\right)  [1+ O(1/z)], \ z \to \infty,\   |\arg z| <  \pi,\eqno(1.8)$$
$$J_\nu(z) = O( z^{\nu} ), \ z \to 0,\eqno(1.9)$$

$$ Y_\nu(z) = \sqrt{2\over \pi z} \sin \left( z- {\pi\over 4} (2\nu+1)\right)  [1+ O(1/z)], \ z \to \infty,\   |\arg z| <  \pi,\eqno(1.10)$$
$$Y_\nu(z) = O\left( |z|^{-|{\rm Re}\nu|} \right), \ z \to 0,\ \nu\neq 0,\eqno(1.11)$$
$$Y_0(z) = O\left(\log(|z|)\right),  \ z \to 0.\eqno(1.12)$$
On the other hand, the modified Bessel functions $I_\nu(z),\ K_\nu(z),\ z,\nu \in \mathbb{C}$ are solutions of the modified Bessel differential equation 
$$  z^2{d^2u\over dz^2}  + z{du\over dz} - (z^2+ \nu^2)u = 0,\eqno(1.13)$$
having the corresponding asymptotic behavior 
$$I_\nu(z) = O\left( |z|^{{\rm Re}\nu} \right), \ z \to 0,\eqno(1.14)$$
$$I_\nu(z) = O\left( {e^z\over \sqrt{2\pi z} } \right), \ z \to \infty,   \ - {\pi\over 2} < \arg z <  {3\pi\over 2},\eqno(1.15)$$
$$K_\nu(z) = O\left( |z|^{-|{\rm Re}\nu|} \right), \ z \to 0,\ \nu\neq 0,\ K_0(z) = O\left(\log(|z|)\right),  \ z \to 0,\eqno(1.16)$$
$$K_\nu(z) = O\left( \sqrt{\pi\over 2z}\  e^{-z} \right), \ z \to \infty, \  | \arg z| <  {3\pi\over 2}.\eqno(1.17)$$
Bessel functions are related by the equalities 

$$Y_{\nu}(z)= {1\over \sin(\pi\nu)}\left[J_\nu(z)\cos(\pi\nu)- J_{-\nu}(z)\right],\eqno(1.18)$$

$$K_\nu(z)= {\pi\over 2\sin(\pi\nu)} \left[ I_{-\nu}(z)- I_\nu(z)\right].\eqno(1.19)$$
Meanwhile, considering the inhomogeneous Bessel equation
$$  z^2{d^2u\over dz^2}  + z{du\over dz} + (z^2- \nu^2)u = z^\mu,\eqno(1.20)$$
we find its solution as the Lommel function $s_{\mu,\nu}(z)$ [2], Ch. 11. Its companion $S_{\mu,\nu}(z)$ is defined by the equality (see [2], Entry 11.9.5)

$$S_{\mu,\nu}(z) = s_{\mu,\nu}(z) + 2^{\mu-1} \Gamma\left( {1\over 2} \left(\mu+\nu+1\right) \right) \Gamma\left( {1\over 2} \left(\mu-\nu+1\right) \right) $$

$$\times\left( \sin \left({1\over 2} \left(\mu-\nu\right)\pi\right) J_\nu(z) - \cos \left({1\over 2} \left(\mu-\nu\right)\pi\right) Y_\nu(z)\right),\eqno(1.21)$$
where $\Gamma(z)$ is Euler's gamma function [2], Ch. 5 and $\mu\pm \nu \neq -1,-2,\dots,$ . It behaves at infinity by virtue of [2], Entry 11.9.9 as follows

$$S_{\mu,\nu}(z)= O\left( z^{\mu-1} \right),\ z \to \infty.\eqno(1.22)$$
In the sequel we will provide existence conditions for discrete transformations (1.1)-(1.6) and establish their inversion formulas for suitable sequences and functions.  To do this,  we will employ integral representations of the Bessel and Lommel functions in the kernels of these operators and classical Fourier series for Lipschitz functions.

\section{Inversion theorems} 

We begin with

{\bf Theorem 1}. {\it   Let a sequence $ \{a_n\}_{n\in \mathbb{N}} $ satisfy the condition  

$$\sum_{n=1}^\infty  |a_n| n  < \infty.\eqno(2.1)$$
Then the discrete transformation $(1.1)$ can be inverted by the formula

$$a_n =  {2\over \pi}  \int_0^\infty   \Phi_n(x) f(x) dx,\ n \in \mathbb{N}_0,\eqno(2.2)$$
where the kernel $\Phi_n(x)$ is defined by 

$$\Phi_n(x) = \int_0^\pi \sin(x\cosh(u)) \sinh(u) \cos(nu) du,\quad x >0,\ n \in \mathbb{N}_0,\eqno(2.3)$$
and  integral  $(2.2)$ converges in the improper sense. }

\begin{proof}  To proceed the proof,  we will appeal to the relatively convergent integral (see [3], Vol. II, Entry 2.12.15.3)

$$\int_0^\infty \sin(t\cosh(u))  {\rm Re} \left[ J_{in}(t) \right] dt =  { \cos(nu) \cosh(\pi n/2)\over \sinh(u)},\ u \in (0,\pi],\  n \in \mathbb{N}_0.\eqno(2.4)$$
In fact,  taking some $T >0$, we have from (1.1)

$$ \int_0^T \sin(t\cosh(u)) \sum_{m=0}^\infty {a_m \over \cosh(\pi m/2)} \ {\rm Re} \left[ J_{im}(t) \right] dt$$

$$ =  \sum_{m=0}^\infty {a_m \over \cosh(\pi m/2)}  \int_0^T \sin(t\cosh(u))  {\rm Re} \left[ J_{im}(t) \right] dt,\eqno(2.5)$$
where the interchange of the order of integration and summation can be justified via  the uniform convergence with respect to $t \in [0,T]$ of the series (1.1).  Indeed, this comes immediately from the definition of the Bessel function and the condition (2.1), because

$$ | J_{in}(t) | =  \left|  \sum_{k=0}^\infty {(-1)^k (t/2)^{2k+in} \over k! \Gamma (k+in+1)} \right| \le  \sum_{k=0}^\infty { (T/2)^{2k} \over k! |\Gamma (k+in+1)|} \le e^T \sqrt{{\sinh(\pi n)\over \pi n}},$$
and therefore

$$ \sum_{m=0}^\infty {|a_m| \over \cosh(\pi m/2)} \ \left| {\rm Re} \left[ J_{im}(t) \right] \right| \le e^T  \sum_{m=0}^\infty  |a_m|\  \sqrt{{2 \tanh(\pi m/2)\over \pi m}} \le e^T ||a||_1 < \infty.\eqno(2.6)$$
Moreover,  taking (2.3),   the right-hand side of (2.2) can be written in the form 

$$\int_0^\infty   \Phi_n(x) f(x) dx = \lim_{T\to \infty} \int_{1/T}^T   \int_0^\pi \sin(x\cosh(u)) \sinh(u) \cos(nu) du $$

$$\times \sum_{m=0}^\infty {a_m \over \cosh(\pi m/2)}  \ {\rm Re} \left[ J_{im}(x) \right] dx $$

$$=  \lim_{T\to \infty}  \sum_{m=0}^\infty {a_m \over \cosh(\pi m/2)}    \int_0^\pi  \int_{1/T}^T  \sin(x\cosh(u)) \sinh(u) \cos(nu) \ {\rm Re} \left[ J_{im}(x) \right] dx du,\eqno(2.7)$$
where the latter interchange of the order of integration and summation is via (2.6).  Now we  employ the integral representation of the kernel in (1.1) (cf. Entry 2.5.54.7 in [3], Vol. I)

$${{\rm Re} \left[ J_{in}(x) \right]\over  \cosh(\pi n/2)} = {2\over \pi} \int_0^\infty \cos(nt) \sin (x\cosh(t)) dt\eqno(2.8)$$
to substitute in (2.7), getting the equality

$$\int_0^\infty   \Phi_n(x) f(x) dx = \lim_{T\to \infty}  {2\over \pi} \sum_{m=0}^\infty a_m   \int_0^\pi  \int_{1/T}^T  \sin(x\cosh(u)) \sinh(u) \cos(nu) $$

$$\times \int_0^\infty \cos(mt) \sin (x\cosh(t)) dt dx du.\eqno(2.9)$$
 Since for some fixed $T >1$,  $ 1/T \le  x \le T$ and sufficiently big $M > T$ we have via integration by parts
 
 $$ \left| \int_M^\infty \cos(nt) \sin (x\cosh(t)) dt\right| = \left| \int_{\cosh(M)}^\infty \cos\left(n \log\left( t + \sqrt {t^2 -1}\right)\right) {\sin (x t)\over \sqrt{ t^2- 1}} dt\right|$$

$$= {1\over x} \left| {\cos(n M) \cos(x \cosh(M)) \over \sinh(M)} - \int_{\cosh(M)}^\infty \cos\left(n \log\left( t + \sqrt {t^2 -1}\right)\right) {\cos (x t) t \over ( t^2- 1)^{3/2}} dt\right.$$

$$\left. + n \int_{\cosh(M)}^\infty \sin \left(n \log\left( t + \sqrt {t^2 -1}\right)\right) {\cos (x t)  \over  t^2- 1} dt\right| \le T \left[ {1 \over \sinh(M)}\right.$$

$$\left. +  \int_{\cosh(M)}^\infty { t \over ( t^2- 1)^{3/2}} dt + n \int_{\cosh(M)}^\infty {1 \over  t^2- 1} dt\right] $$

$$= T \left[ {2 \over \sinh(M)} + {n\over 2} \log\left({\cosh(M) +1\over \cosh(M)- 1}\right) \right] \to 0,\quad M \to \infty.$$
Therefore integral (2.8) converges uniformly with respect to $ x \in [1/T, T].$ Consequently, we change the order of integration and calculating an elementary integral, we  obtain

$$\int_0^\infty   \Phi_n(x) f(x) dx = \lim_{T\to \infty}  {1\over \pi} \sum_{m=0}^\infty a_m   \int_0^\pi  \sinh(u) \cos(nu) $$

$$\times \int_0^\infty \cos(mt) \left[ {\sin (T\left(\cosh(t) - \cosh(u)\right))\over \cosh(t) - \cosh(u)}-  {\sin (T\left(\cosh(t) + \cosh(u)\right))\over \cosh(t) + \cosh(u)}\right. $$

$$\left. - {\sin (1/T\left(\cosh(t) - \cosh(u)\right))\over \cosh(t) - \cosh(u)} +  {\sin (1/T\left(\cosh(t) + \cosh(u)\right))\over \cosh(t) + \cosh(u)}\right] dt  du .\eqno(2.10)$$
In the meantime, since via (2.1) 

$$\left|\sum_{m=0}^\infty a_m   \int_0^\pi  \sinh(u) \cos(nu) \right. $$

$$\left.\times \int_0^\infty \cos(mt) \left[   {\sin (1/T\left(\cosh(t) - \cosh(u)\right))\over \cosh(t) - \cosh(u)} - {\sin (1/T\left(\cosh(t) + \cosh(u)\right))\over \cosh(t) + \cosh(u)}\right] dt  du \right|$$

$$=  \left|\sum_{m=0}^\infty a_m   \int_0^\pi  \sinh(u) \cos(nu) \right. $$

$$\left.\times \int_0^{\sqrt T} \cos(mt) \left[   {\sin (1/T\left(\cosh(t) - \cosh(u)\right))\over \cosh(t) - \cosh(u)} - {\sin (1/T\left(\cosh(t) + \cosh(u)\right))\over \cosh(t) + \cosh(u)}\right] dt  du \right.$$

$$+ \left.\sum_{m=0}^\infty a_m   \int_0^\pi  \sinh(u) \cos(nu) \right. $$

$$\left.\times \int_{\sqrt T}^\infty \cos(mt) \left[   {\sin (1/T\left(\cosh(t) - \cosh(u)\right))\over \cosh(t) - \cosh(u)} - {\sin (1/T\left(\cosh(t) + \cosh(u)\right))\over \cosh(t) + \cosh(u)}\right] dt  du \right|$$

$$\le  [\cosh(\pi)-1] \sum_{m=0}^\infty | a_m | \left[ {2\over \sqrt T} +  \int_{\sqrt T}^\infty  \left[   {1\over \cosh(t) - \cosh(\pi)} + {1\over \cosh(t)}\right] dt \right] \to 0,\ T \to \infty,$$
our goal will be to justify the existence of the limit and pass to it under series sign in the equality 

$$\int_0^\infty   \Phi_n(x) f(x) dx = \lim_{T\to \infty}  {1\over \pi} \sum_{m=0}^\infty a_m   \int_0^\pi  \sinh(u) \cos(nu) $$

$$\times \int_0^\infty \cos(mt) \left[ {\sin (T\left(\cosh(t) - \cosh(u)\right))\over \cosh(t) - \cosh(u)}-  {\sin (T\left(\cosh(t) + \cosh(u)\right))\over \cosh(t) + \cosh(u)}\right] dt du.\eqno(2.11)$$
In fact, since

$$ \sum_{m=0}^\infty |a_m |  \int_0^\pi  \sinh(u) |\cos(nu)|  \int_0^\infty\left| \cos(mt)  {\sin (T\left(\cosh(t) + \cosh(u)\right))\over \cosh(t) + \cosh(u)}\right| dt du $$

$$\le  \sum_{m=0}^\infty |a_m |   \int_0^\pi  \sinh(u)   du \int_0^\infty {dt\over \cosh(t)}=  {\pi\over 2}  [\cosh(\pi)-1]  \sum_{m=0}^\infty |a_m | < \infty, $$
we have 

$$ \lim_{T\to \infty} \sum_{m=0}^\infty a_m   \int_0^\pi  \sinh(u) \cos(nu)  \int_0^\infty \cos(mt) \  {\sin (T\left(\cosh(t) + \cosh(u)\right))\over \cosh(t) + \cosh(u)}dt du$$

$$=  \int_0^\pi  \sinh(u) \cos(nu)   \lim_{T\to \infty} \int_{1+\cosh(u)}^\infty  \sum_{m=0}^\infty a_m  \cos\left(m \log\left( t -\cosh(u) + ( (t-\cosh(u))^2 -1)^{1/2}\right)\right)$$

$$\times   {\sin (T t)\over t  ((t-\cosh(u))^2 -1)^{1/2}}dt du = 0,\eqno(2.12)$$
owing to the Riemann-Lebesgue lemma.  Further,  we write 
 
 $${1\over \pi} \sum_{m=0}^\infty a_m   \int_0^\pi  \sinh(u) \cos(nu) \int_0^\infty \cos(mt) \  {\sin (T\left(\cosh(t) - \cosh(u)\right))\over \cosh(t) - \cosh(u)} dt du $$

 $$= {1\over \pi} \sum_{m=0}^\infty a_m   \int_0^{\cosh(\pi)-1} \cos \left(n \log\left( u+ 1+( (u+1)^2 -1)^{1/2}\right) \right) $$
 
 $$\times \int_{-u}^\infty  \cos \left(m \log\left( t +u +1+ ( (t+u+1)^2 -1)^{1/2}\right) \right) {\sin (T t)\over t ((t+u+1)^2-1)^{1/2}} dt du $$

 $$= {1\over \pi} \sum_{m=0}^\infty a_m   \int_0^{\cosh(\pi)-1} \cos \left(n \log\left( u+ 1+( (u+1)^2 -1)^{1/2}\right) \right) $$
 
 $$\times \int_{-u}^\infty  \left[{\cos \left(m \log\left( t +u +1+ ( (t+u+1)^2 -1)^{1/2}\right) \right)\over  ((t+u+1)^2-1)^{1/2}} \right.$$
 
 $$\left. - {\cos \left(m \log\left( u +1+ ( (u+1)^2 -1)^{1/2} \right)\right)\over  ((u+1)^2-1)^{1/2}} \right]   {\sin (T t)\over t} dt du $$
 
 $$+ {1\over \pi} \sum_{m=0}^\infty a_m   \int_0^{\cosh(\pi)-1} \cos \left(n \log\left( u+ 1+( (u+1)^2 -1)^{1/2}\right) \right) $$
 
 $$\times   {\cos \left(m \log\left( u +1+ ( (u+1)^2 -1)^{1/2} \right)\right)\over   ((u+1)^2-1)^{1/2}} \int_{-u}^\infty {\sin (T t)\over t} dt du.\eqno(2.13) $$
Then, in turn, 

$${1\over \pi} \sum_{m=0}^\infty a_m   \int_0^{\cosh(\pi)-1} \cos \left(n \log\left( u+ 1+( (u+1)^2 -1)^{1/2}\right) \right) $$
 
 $$\times   {\cos \left(m \log\left( u +1+ ( (u+1)^2 -1)^{1/2} \right)\right)\over   ((u+1)^2-1)^{1/2}} \int_{-u}^\infty {\sin (T t)\over t} dt du $$

 $$= {1\over \pi} \sum_{m=0}^\infty a_m  \left( \int_0^{1/\sqrt T} + \int_{1/\sqrt T}^{\cosh(\pi)-1} \right) \cos \left(n \log\left( u+ 1+( (u+1)^2 -1)^{1/2}\right) \right) $$
 
 $$\times   {\cos \left(m \log\left( u +1+ ( (u+1)^2 -1)^{1/2} \right)\right)\over   ((u+1)^2-1)^{1/2}} \int_{-u}^\infty {\sin (T t)\over t} dt du $$
 and, employing the second mean value theorem,
 
 $$ {1\over \pi} \sum_{m=0}^\infty a_m  \int_{1/\sqrt T}^{\cosh(\pi)-1}  \cos \left(n \log\left( u+ 1+( (u+1)^2 -1)^{1/2}\right) \right) $$
 
 $$\times   {\cos \left(m \log\left( u +1+ ( (u+1)^2 -1)^{1/2} \right)\right)\over   ((u+1)^2-1)^{1/2}} \int_{-u}^\infty {\sin (T t)\over t} dt du $$

 $$=  {1\over \pi} \sum_{m=0}^\infty a_m  \int_{1/\sqrt T}^{\cosh(\pi)-1}  \cos \left(n \log\left( u+ 1+( (u+1)^2 -1)^{1/2}\right) \right) $$
 
 $$\times   {\cos \left(m \log\left( u +1+ ( (u+1)^2 -1)^{1/2} \right)\right)\over   ((u+1)^2-1)^{1/2}} \left( \pi- \int_{u T}^\infty {\sin (t)\over t} dt \right) du $$

 $$=   \sum_{m=0}^\infty a_m  \int_{1/\sqrt T}^{\cosh(\pi)-1}  \cos \left(n \log\left( u+ 1+( (u+1)^2 -1)^{1/2}\right) \right) $$
 
 $$\times   {\cos \left(m \log\left( u +1+ ( (u+1)^2 -1)^{1/2} \right)\right)\over   ((u+1)^2-1)^{1/2}} du $$

 $$-  {1\over T\pi} \sum_{m=0}^\infty a_m  \int_{1/\sqrt T}^{\cosh(\pi)-1}  \cos \left(n \log\left( u+ 1+( (u+1)^2 -1)^{1/2}\right) \right) $$
 
 $$\times   {\cos \left(m \log\left( u +1+ ( (u+1)^2 -1)^{1/2} \right)\right)\over u  ((u+1)^2-1)^{1/2}} \int_{u T}^{T^\prime} \sin (t) dt du  $$

 $$\to \sum_{m=0}^\infty a_m  \int_{0}^{\pi}  \cos \left(n u \right) \cos \left(m u \right) du = {\pi\over 2} a_n,\quad T \to \infty\eqno(2.14)$$
because

$$ {1\over T\pi} \left| \sum_{m=0}^\infty a_m  \int_{1/\sqrt T}^{\cosh(\pi)-1}  \cos \left(n \log\left( u+ 1+( (u+1)^2 -1)^{1/2}\right) \right)\right. $$
 
 $$\left. \times   {\cos \left(m \log\left( u +1+ ( (u+1)^2 -1)^{1/2} \right)\right)\over u  ((u+1)^2-1)^{1/2}} \int_{u T}^{T^\prime} \sin (t) dt du \right| $$

$$\le {A \over (T (1+2\sqrt T))^{1/2}}  \sum_{m=0}^\infty  |a_m|   \int_{1/\sqrt T}^{\cosh(\pi)-1}  {du\over u}$$

$$ = {A \log(T^{1/2} (\cosh(\pi)-1)) \over (T (1+2\sqrt T))^{1/2}}  \sum_{m=0}^\infty  |a_m| \to 0,\ T \to \infty,$$
where $ A >0$ is an absolute constant.  Hence, returning to (2.13) we need to establish the equality 

$$ {1\over \pi} \lim_{T\to \infty} \sum_{m=0}^\infty a_m   \int_0^{\cosh(\pi)-1} \cos \left(n \log\left( u+ 1+( (u+1)^2 -1)^{1/2}\right) \right) $$
 
 $$\times \int_{-u}^\infty  \left[{\cos \left(m \log\left( t +u +1+ ( (t+u+1)^2 -1)^{1/2}\right) \right)\over  ((t+u+1)^2-1)^{1/2}} \right.$$
 
 $$\left. - {\cos \left(m \log\left( u +1+ ( (u+1)^2 -1)^{1/2} \right)\right)\over  ((u+1)^2-1)^{1/2}} \right]   {\sin (T t)\over t} dt du = 0.\eqno(2.15)$$
In fact, denoting the expression under the limit sign in (2.15) by $I(T)$ and fixing a small positive $\delta$, we get 

$$I(T)=  {1\over \pi}  \sum_{m=0}^\infty a_m   \int_0^{\cosh(\pi)-1} \cos \left(n \log\left( u+ 1+( (u+1)^2 -1)^{1/2}\right) \right) $$
 
 $$\times \int_{-u}^\delta  \left[{\cos \left(m \log\left( t +u +1+ ( (t+u+1)^2 -1)^{1/2}\right) \right)\over  ((t+u+1)^2-1)^{1/2}} \right.$$
 
 $$\left. - {\cos \left(m \log\left( u +1+ ( (u+1)^2 -1)^{1/2} \right)\right)\over  ((u+1)^2-1)^{1/2}} \right]   {\sin (T t)\over t} dt du $$

$$+  {1\over \pi}  \sum_{m=0}^\infty a_m   \int_0^{\cosh(\pi)-1} \cos \left(n \log\left( u+ 1+( (u+1)^2 -1)^{1/2}\right) \right) $$
 
 $$\times \int_\delta^\infty  \left[{\cos \left(m \log\left( t +u +1+ ( (t+u+1)^2 -1)^{1/2}\right) \right)\over  ((t+u+1)^2-1)^{1/2}} \right.$$
 
 $$\left. - {\cos \left(m \log\left( u +1+ ( (u+1)^2 -1)^{1/2} \right)\right)\over  ((u+1)^2-1)^{1/2}} \right]   {\sin (T t)\over t} dt du= I_1(T)+ I_2(T).$$
Hence,  
 
$$\lim_{T\to \infty} I_2(T) =  {1\over \pi}  \lim_{T\to \infty} \sum_{m=0}^\infty a_m   \int_0^{\cosh(\pi)-1} \cos \left(n \log\left( u+ 1+( (u+1)^2 -1)^{1/2}\right) \right) $$
 
 $$\times \int_\delta^\infty  {\cos \left(m \log\left( t +u +1+ ( (t+u+1)^2 -1)^{1/2}\right) \right)\over  ((t+u+1)^2-1)^{1/2}}  {\sin (T t)\over t} dt du $$
 
 $$ - {1\over \pi}  \lim_{T\to \infty} \sum_{m=0}^\infty a_m   \int_0^{\cosh(\pi)-1} \cos \left(n \log\left( u+ 1+( (u+1)^2 -1)^{1/2}\right) \right) $$
 
 $$\times {\cos \left(m \log\left( u +1+ ( (u+1)^2 -1)^{1/2} \right)\right)\over  ((u+1)^2-1)^{1/2}} du  \int_{\delta T}^\infty {\sin (t)\over t} dt =0$$
 via the Riemann-Lebesgue lemma and the remainder of the convergent integral.  Concerning the expression $I_1(T)$, we appeal to the Leibniz differentiation formula under the integral sign
 to write
 
 $$ I_1(T)=  {a_0\over \pi}  \int_0^{\cosh(\pi)-1} \cos \left(n \log\left( u+ 1+( (u+1)^2 -1)^{1/2}\right) \right) $$
 
 $$\times \int_{-u}^\delta  \left[{1\over  ((t+u+1)^2-1)^{1/2}} - {1\over  ((u+1)^2-1)^{1/2}} \right]   {\sin (T t)\over t} dt du $$

$$+ {1\over \pi}  \sum_{m=1}^\infty {a_m \over m}   \int_0^{\cosh(\pi)-1} \cos \left(n \log\left( u+ 1+( (u+1)^2 -1)^{1/2}\right) \right) $$
 
 $$\times {d\over du} \int_{-u}^\delta  \left[\sin \left(m \log\left( t +u +1+ ( (t+u+1)^2 -1)^{1/2}\right) \right)\right.$$
 
 $$\left. - \sin \left(m \log\left( u +1+ ( (u+1)^2 -1)^{1/2} \right)\right) \right]   {\sin (T t)\over t} dt du $$

 $$+ {1\over \pi}  \sum_{m=0}^\infty a_m   \int_0^{\cosh(\pi)-1} \cos \left(n \log\left( u+ 1+( (u+1)^2 -1)^{1/2}\right) \right) $$
 
 $$\times \sin \left(m \log\left( u +1+ ( (u+1)^2 -1)^{1/2} \right)\right)   {\sin (T u)\over u}  du$$
 
 $$=N_1(T)+ N_2(T)+ N_3(T).\eqno(2.16) $$
 Then since
 
 $$\int_0^{\cosh(\pi)-1} \Bigg|\cos \left(n \log\left( u+ 1+( (u+1)^2 -1)^{1/2}\right) \right) \Bigg.$$
 
 $$\left.\times \int_{-u}^\delta  \left[{1\over  ((t+u)(t+u+2))^{1/2}} - {1\over  (u(u+2))^{1/2}} \right]   {\sin (T t)\over t}\right| dt du $$
 
 $$\le  \int_0^{\cosh(\pi)-1} {1\over (u(u+2))^{3/4}}  \int_{-u}^\delta  {t+2(u+1) \over ((t+u)(t+u+2))^{3/4} } dt du  < \infty,$$
it yields $ \lim_{T\to \infty} N_1(T) =0 $ via the Riemann-Lebesgue lemma.  Analogously, since 
 
 $$\cos \left(n \log\left( u+ 1+( (u+1)^2 -1)^{1/2}\right) \right) $$
 
 $$\times \sin \left(m \log\left( u +1+ ( (u+1)^2 -1)^{1/2} \right)\right)   {1\over u} \in L_1(0, \cosh(\pi)-1),$$
 we get, taking into account (2.1), $ \lim_{T\to \infty} N_3(T) =0 $.  Finally, the middle term in (2.16) can be treated via integration by parts. In fact, we find

 $$N_2(T)=  {(-1)^n \over \pi}  \sum_{m=1}^\infty {a_m \over m}   \int_{1-\cosh(\pi)}^\delta \sin \left(m \log\left( t +\cosh(\pi) + ( (t+ \cosh(\pi))^2 -1)^{1/2}\right) \right) {\sin (T t)\over t} dt $$
 
 $$- {1\over \pi}  \sum_{m=1}^\infty {a_m \over m}  \int_{0}^\delta  \sin \left(m \log\left( t +1+ ( (t+1)^2 -1)^{1/2}\right) \right) {\sin (T t)\over t} dt  $$

$$ + {2n \over \pi}  \sum_{m=1}^\infty {a_m \over m}   \int_0^{\cosh(\pi)-1} {\sin \left(n \log\left( u+ 1+( (u+1)^2 -1)^{1/2}\right) \right) \over ( u(u+2))^{1/2}}$$
 
 $$\times  \int_{-u}^\delta  \sin \left({m\over 2}  \log\left( {t +u +1+ ( (t+u+1)^2 -1)^{1/2}\over  u +1+ ( (u+1)^2 -1)^{1/2}} \right) \right)$$
 
 $$\times \cos \left({m\over 2}  \log\left( (t +u +1+ ( (t+u+1)^2 -1)^{1/2} ) ( u +1+ ( (u+1)^2 -1)^{1/2} ) \right) \right)  {\sin (T t)\over t} dt du, $$
 and we see that   $ \lim_{T\to \infty} N_2(T) =0 $ by the same reasons. Thus, combining with (2.12), (2.13), (2.14), we return to (2.11) to establish the inversion formula ( 2.2), completing the proof of Theorem 1. 
 
 \end{proof}
 
 The discrete transformation (1.2) can be inverted by the following theorem.

{\bf Theorem 2}.   {\it Let $f$ be a complex-valued function on $\mathbb{R}_+$ which is represented by the integral 

$$f(x) =   \int_{0}^\pi \sin (x \cosh(u))  \varphi(u)  du,\quad x >0,\eqno(2.17)$$ 
where $ \varphi(u) = \psi(u)\sinh(u),\ \psi(0)=0$ and $\psi$ is a continuously differentiable  even   $2\pi$-periodic function.  Then for all $x >0$ the following inversion formula for  transformation $(1.2)$  holds

$$ f(x)  = {2\over x\pi}  \sin\left({x\over 2} (\cosh(\pi)-1)\right)\sin\left({x\over 2} (\cosh(\pi)+1)\right)+  {2\over \pi} \sum_{n=0}^\infty   \Phi_n (x) a_n,\eqno(2.18)$$
where $\Phi_n$ is defined by $(2.3)$.}

\begin{proof}    Plugging the right-hand side of the representation (2.17) in (1.2), we change the order of integration,  employ  (2.4)  and the definition of $\varphi$ to obtain 

$$a_n = \int_{0}^\pi \psi (u) \cos(nu) du.\eqno(2.19)$$
The interchange of the order of integration can be justified in the following way.  We have

$$  \int_0^\infty   {\rm Re} \left[ J_{in}(x) \right] f(x) dx = \lim_{T\to \infty} \int_0^T  {\rm Re} \left[ J_{in}(x) \right] f(x) dx $$

$$=  \lim_{T\to \infty} \int_{0}^\pi  \varphi(u) \int_0^T  {\rm Re} \left[ J_{in}(x) \right] \sin (x \cosh(u))  dx  du$$

$$\cosh\left({\pi n\over 2}\right) \int_{0}^\pi \psi (u) \cos(nu) du -  \lim_{T\to \infty} \int_{0}^\pi  \varphi(u) \int_T^\infty  {\rm Re} \left[ J_{in}(x) \right] \sin (x \cosh(u))  dx  du,$$
where the interchange is allowed owing to the continuity of the integrand on the rectangle $[0,\pi] \times [0,T]$.  To show that the latter limit is zero we appeal to the asymptotic behavior of the Bessel function (1.8) to find for each $n \in \mathbb{N}$

$$ {\rm Re} \left[ J_{in}(x) \right] = \cosh\left({\pi n\over 2}\right) \sqrt{2\over \pi x} \cos \left( x- {\pi\over 4} \right)  [1+ O(1/x)],\quad  x \to \infty.$$
Hence for sufficiently big $T > 0$ and by virtue of second mean value theorem we get 
$$\int_{0}^\pi  \varphi(u) \int_T^\infty  {\rm Re} \left[ J_{in}(x) \right] \sin (x \cosh(u))  dx  du$$

$$ = \cosh\left({\pi n\over 2}\right) \sqrt{1\over 2\pi }  \int_{0}^\pi  \varphi(u) \int_T^\infty   \sin \left(x (\cosh(u)-1) - {\pi\over 4} \right) {dx du\over \sqrt x}  + O\left({1\over \sqrt T}\right)$$

$$= O\left( {1\over \sqrt T}   \int_{0}^\pi  {\varphi(u)\over \cosh(u)-1} du \right) \to 0,\ T \to \infty$$
since the latter integral is convergent under conditions of the theorem.  Thus, returning to (2.19),  we substitute $a_n$ and $\Phi_n$ by (2.3) into the partial sum of the series (2.18) $S_N(x) $, and  it becomes 

$$ S_N(x)  = {2\over \pi} \sum_{n=0}^N   \int_{0}^\pi \sin(x \cosh(t))  \sinh(t) \cos(nt) dt \int_{0}^\pi   \psi(u) \cos(nu) du.\eqno(2.20)$$
Hence,  calculating the sum via the known identity 
$$ \sum_{n=0}^N  \cos(nt) \cos(nu) = {1\over 4} \left[ 2+  {\sin \left((2N+1) (u-t)/2 \right)\over \sin( (u-t) /2)}  +  {\sin \left((2N+1) (u+t)/2 \right)\over \sin( (u+t) /2)} \right],$$
 we obtain from (2.20) 
 
 $$ S_N(x)  = {2\over x\pi}  \sin\left({x\over 2} (1-\cosh(\pi))\right)\sin\left({x\over 2} (\cosh(\pi)+1)\right) +  {1\over 2\pi} \int_{0}^\pi \sin(x \cosh(t)) \sinh(t) $$
 
 $$\times  \int_{-\pi}^\pi    \psi(u)  {\sin \left((2N+1) (u-t)/2 \right)\over \sin( (u-t) /2)} du dt.\eqno(2.21)$$
  Since $\psi$ is $2\pi$-periodic, we treat  the latter integral with respect to $u$ as follows 

$$  \int_{-\pi}^{\pi}  \psi(u) \ {\sin \left((2N+1) (u-t)/2 \right)\over \sin( (u-t) /2)}  du $$

$$=  \int_{ t-\pi}^{t+ \pi} \psi(u) \  {\sin \left((2N+1) (u-t)/2 \right)\over \sin( (u-t) /2)}  du $$

$$=  \int_{ -\pi}^{\pi}   \psi(u+t) \  {\sin \left((2N+1) u/2 \right)\over \sin( u /2)}  du. $$
Moreover,

$$ {1\over 2\pi} \int_{ -\pi}^{\pi}  \psi(u+t) \  {\sin \left((2N+1) u/2 \right)\over \sin( u /2)}  du -  \psi(t) $$

$$=  {1\over 2\pi} \int_{ -\pi}^{\pi}  \left[ \psi(u+t)- \psi(t) \right]  \  {\sin \left((2N+1) u/2 \right)\over \sin( u /2)}  du.$$
When  $u+t > \pi$ or  $u+t < -\pi$ then we interpret  the value  $\psi(u+t)$ by  formulas

$$\psi(u+t)- \psi(t)= \psi(u+t-2\pi)- \psi(t - 2\pi),$$ 

$$\psi(u+t)- \psi(t) = \psi(u+t+ 2\pi)- \psi(t +2\pi),$$ 
respectively.     Then   since $\psi \in C^1[-\pi,\pi]$,  it satisfies the Lipschitz condition on $[-\pi, \pi]$

$$\left| \psi(u) - \psi(v)\right| \le C |u-v|, \quad  \forall \  u, v \in  [-\pi, \pi],\eqno(2.22)$$
where $C >0$ is an absolute constant.  Hence  we have the uniform estimate for any $t \in [-\pi,\pi]$

$${\left|  \psi(u+t)- \psi(t) \right| \over | \sin( u /2) |}  \le C \left| {u\over \sin( u /2)} \right|.$$
Therefore,  owing to the Riemann-Lebesgue lemma

$$\lim_{N\to \infty } {1\over 2\pi} \int_{ -\pi}^{\pi}  \left[ \psi(u+t)  - \psi(t) \right]  \  {\sin \left((2N+1) u/2 \right)\over \sin( u /2)}  du =  0\eqno(2.23)$$
for all $ t\in [-\pi,\pi].$    Besides, returning to (2.21), we estimate the iterated integral 
$$ \int_{0}^\pi \left|\sin(x \cosh(t)) \right| \sinh(t)  \int_{ -\pi}^{\pi} \left| \left[ \psi(u+t)  - \psi(t)  \right]\right.$$

$$\left.\times   {\sin \left((2N+1) u/2 \right)\over \sin( u /2)}  \right| du dt \le  C [\cosh(\pi)-1]  \int_{ -\pi}^{\pi}   \left| {u\over \sin( u /2)} \right| du < \infty.$$
 Consequently, via  the dominated convergence theorem it is possible to pass to the limit when $N \to \infty$ under the  integral sign, and recalling (2.23), we derive

$$  \lim_{N \to \infty}   {1\over 2\pi}  \int_{0}^\pi \sin (x \cosh(t))  \sinh(t)   \int_{ -\pi}^{\pi}  \left[ \psi(u+t)  - \psi(t)  \right] $$

$$\times  \  {\sin \left((2N+1) u/2 \right)\over \sin( u /2)}  du dt =  {1\over 2 \pi}  \int_{0}^\pi \sin (x \cosh(t))  \sinh(t)   $$

$$ \times \lim_{N \to \infty}  \int_{ -\pi}^{\pi}  \left[ \psi(u+t)-   \psi(t)  \right]  \  {\sin \left((2N+1) u/2 \right)\over \sin( u /2)}  du dt = 0.$$
Hence, combining with (2.21),  we obtain  by virtue of  the definition of $\varphi$ and $f$

$$ \lim_{N \to \infty}  S_N(x) =  {2\over x\pi}  \sin\left({x\over 2} (1-\cosh(\pi))\right)\sin\left({x\over 2} (\cosh(\pi)+1)\right) +  \int_{0}^\pi  \sin(x \cosh(t))  \varphi (t) dt$$

$$ = f(x) + {2\over x\pi}  \sin\left({x\over 2} (1-\cosh(\pi))\right)\sin\left({x\over 2} (\cosh(\pi)+1)\right),$$
where the integral (2.17) converges since $\varphi \in C[0,\pi]$.  Thus we established  (2.18), completing the proof of Theorem 2.
 
\end{proof}

The same scheme can be applied to invert discrete $Im$-transformations (1.3), (1.4).  It involves the following analogs of the integrals (2.4), (2.8) (cf. [3], Vol. II, Entry 2.12.15.3, Vol. I, Entry 2.5.54.7)

$$\int_0^\infty \cos(t\cosh(u))  {\rm Im} \left[ J_{in}(t) \right] dt =  - { \cos(nu) \sinh(\pi n/2)\over \sinh(u)},\ u \in (0,\pi],\  n \in \mathbb{N},\eqno(2.24)$$

$${{\rm Im} \left[ J_{in}(x) \right]\over  \sinh(\pi n/2)} = - {2\over \pi} \int_0^\infty \cos(nt) \cos (x\cosh(t)) dt.\eqno(2.25)$$
We will formulate the corresponding  theorems, leaving the proofs to the interested reader. 

{\bf Theorem 3}. {\it   Let a sequence $ \{a_n\}_{n\in \mathbb{N}} $ satisfy the condition  $(2.1)$.  Then  the discrete transformation $(1.3)$ can be inverted by the formula

$$a_n =  - {2\over \pi}  \int_0^\infty   \Psi_n(x) f(x) dx,\ n \in \mathbb{N}_0,\eqno(2.26)$$
where the kernel $\Psi_n(x)$ is defined by 

$$\Psi_n(x) = \int_0^\pi \cos(x\cosh(u)) \sinh(u) \cos(nu) du,\quad x >0,\ n \in \mathbb{N}_0,\eqno(2.27)$$
and  integral  $(2.26)$ converges in the improper sense. }

{\bf Theorem 4}.   {\it Let $f$ be a complex-valued function on $\mathbb{R}_+$ which is represented by the integral 

$$f(x) =   \int_{0}^\pi \cos (x \cosh(u))  \varphi(u)  du,\quad x >0,\eqno(2.28)$$ 
where $ \varphi(u) = \psi(u)\sinh(u),\  \psi(0)=0$ and $\psi$ is a continuously differentiable even   $2\pi$-periodic function.  Then for all $x >0$ the following inversion formula for  transformation $(1.4)$  holds

$$ f(x)  = {2\over x\pi}  \sin\left({x\over 2} (\cosh(\pi)-1)\right)\cos\left({x\over 2} (\cosh(\pi)+1)\right)+  {2\over \pi} \sum_{n=0}^\infty   \Psi_n (x) a_n,\eqno(2.29)$$
where $\Psi_n$ is defined by $(2.27)$.}

In order to establish the inversion formula for the transformation (1.5), involving the Lommel function (1.21) in the kernel, we will set up the following lemma.

{\bf Lemma 1}. {\it  Let  ${\rm Re}\mu < 0, u \in (0,\pi],\ n \in \mathbb{N}$. Then the following formula takes place

$$\int_0^\infty \sin\left(x\cosh(u)- {\pi\mu\over 2}\right) S_{\mu,in} (x) dx $$

$$= {2^{\mu} \pi^2 \sin(nu) \over \sinh(u)\sinh(\pi n) \Gamma((1-\mu-in)/2) \Gamma((1-\mu+in)/2)},\eqno(2.30)$$
and integral $(2.30)$ converges absolutely.}

\begin{proof}  Taking the Mellin-Barnes representation for the Lommel function $S_{\mu,in} (x)$ (see [3], Vol. III, Entry 8.4.27.3) with the use of the reflection formula for the gamma function, we find

$$ S_{\mu,in} (2 x) = {2^{\mu-1} \over 4i \Gamma((1-\mu-in)/2) \Gamma((1-\mu+in)/2)}$$

$$\times  \int_{\gamma-i\infty}^{\gamma+i\infty} \Gamma\left({s+in\over 2}\right) \Gamma\left({s-in\over 2}\right) {x^{-s}  \over \cos (\pi (s+\mu)/2)} ds, \eqno(2.31)$$
where $x > 0,\ -1- {\rm Re} \mu, 0 < \gamma < 1- {\rm Re} \mu$.  Then,  shifting the contour in the integral (2.31) to the left and to the right within this vertical strip and appealing  to the Stirling asymptotic formula for the gamma function (see [2], Entry 5.11.9),  one can guarantee the convergence of the following iterated integral

$$\left(\int_0^1+ \int_1^\infty\right) \left[\  \left|\cos\left({\pi\mu\over 2}\right)\sin(2x\cosh(u)) \right| +  \left|\sin\left({\pi\mu\over 2}\right)\cos(2x\cosh(u)) \right| \right]$$

$$\times \int_{\gamma-i\infty}^{\gamma+i\infty} \left|\Gamma\left({s+in\over 2}\right) \Gamma\left({s-in\over 2}\right) {x^{-s}  \over \cos (\pi (s+\mu)/2)} ds\right| dx$$

$$\le \int_0^1\left[\  \left|\cos\left({\pi\mu\over 2}\right)\right| +  \left|\sin\left({\pi\mu\over 2}\right) \right| \right] {dx\over x^\gamma}$$

$$\times \int_{\gamma-i\infty}^{\gamma+i\infty} \left|\Gamma\left({s+in\over 2}\right) \Gamma\left({s-in\over 2}\right) {ds \over \cos (\pi (s+\mu)/2)} \right| $$

$$+  \int_1^\infty \left[\  \left|\cos\left({\pi\mu\over 2}\right)\right| +  \left|\sin\left({\pi\mu\over 2}\right) \right| \right] {dx\over x^\gamma}$$

$$\times \int_{\gamma-i\infty}^{\gamma+i\infty} \left|\Gamma\left({s+in\over 2}\right) \Gamma\left({s-in\over 2}\right) {ds \over \cos (\pi (s+\mu)/2)} \right|  < \infty,$$
where we choose $\gamma \in (\max( -1- {\rm Re} \mu, 0 ), 1)$ for the integration over $(0,1)$  and $\gamma \in (1,  1- {\rm Re} \mu )$ for the integration over $(1, \infty).$
Therefore, making use  the integral (2.31), we deduce from  (2.30) after the interchange of  the order of integration by Fubini's theorem and involving  Entries 8.4.5.1, 8.4.5.2 in [3], Vol. III 

$$\int_0^\infty \sin\left(2x\cosh(u)- {\pi\mu\over 2}\right)  S_{\mu,in} (2x) dx $$

$$=  {2^{\mu-1} \over 4i \Gamma((1-\mu-in)/2) \Gamma((1-\mu+in)/2)}$$

$$\times  \int_{\gamma-i\infty}^{\gamma+i\infty} \Gamma\left({s+in\over 2}\right) \Gamma\left({s-in\over 2}\right) {1  \over \cos (\pi (s+\mu)/2)} $$

$$\times  \int_0^\infty \sin\left(2x\cosh(u)- {\pi\mu\over 2}\right) x^{-s}  dx ds$$

$$=  {2^{\mu-2} \sqrt \pi\over 4i \Gamma((1-\mu-in)/2) \Gamma((1-\mu+in)/2)}$$

$$\times  \int_{\gamma-i\infty}^{\gamma+i\infty} \Gamma\left({s+in\over 2}\right) \Gamma\left({s-in\over 2}\right) { (\cosh(u))^{s-1} \over \cos (\pi (s+\mu)/2)} $$

$$\times   \left[ \cos\left({\pi\mu\over 2}\right) {\Gamma(1- s/2)\over \Gamma((1+s)/2)}   -  \sin\left({\pi\mu\over 2}\right)   {\Gamma((1- s)/2)\over \Gamma(s/2)} \right]  x^{-s} ds$$

$$=  {2^{\mu-2} \over 4i \Gamma((1-\mu-in)/2) \Gamma((1-\mu+in)/2)}$$

$$\times  \int_{\gamma-i\infty}^{\gamma+i\infty} \Gamma(1-s) \Gamma\left({s+in\over 2}\right) \Gamma\left({s-in\over 2}\right)  2^{s} (\cosh(u))^{s-1} ds.\eqno(2.32)$$
But the latter integral can be treated via the Parseval equality for the Mellin transform [4] and Entries 8.4.3.1, 8.4.23.1 in [2], Vol. III, 2.16.6.1 in [2], Vol. II. Thus we obtain 

$${2^{\mu-2} \over 4i \Gamma((1-\mu-in)/2) \Gamma((1-\mu+in)/2)}$$

$$\times  \int_{\gamma-i\infty}^{\gamma+i\infty} \Gamma(1-s) \Gamma\left({s+in\over 2}\right) \Gamma\left({s-in\over 2}\right)  2^{s} (\cosh(u))^{s-1} ds$$

$$ = {2^{\mu-1} \pi \over  \Gamma((1-\mu-in)/2) \Gamma((1-\mu+in)/2)} \int_0^\infty e^{- x\cosh(u)} K_{in} (x) dx$$

$$ = {2^{\mu-1} \pi^2 \sin(nu) \over \sinh(u)\sinh(\pi n) \Gamma((1-\mu-in)/2) \Gamma((1-\mu+in)/2)} .$$
Thus, combining with (2.32), we arrive at (2.30), completing the proof of Lemma 1.

\end{proof}

{\bf Corollary 1}. {\it  Let $x >0, - 5/4 < {\rm Re}\mu < 3/4,\ n \in \mathbb{N}$. The following inequality holds valid

$$\left| S_{\mu,in} (x)\right| \le  { C\ x^{-1/4} \over [\sinh(\pi n)]^{1/2}   |\Gamma((1-\mu-in)/2) \Gamma((1-\mu+in)/2)|},\eqno(2.33)$$ 
where $C > 0$ is an absolute constant.}

\begin{proof} To prove (2.33), we appeal to the Lebedev inequality for the modified Bessel function (cf. [4], p.219)

$$|K_{in}(x)| \le A\  {x^{-1/4} \over [\sinh(\pi n)]^{1/2} },\  x >0, n \in \mathbb{N},\eqno(2.34)$$
where $A > 0$ is an absolute constant. Then from (2.31) and the Parseval equality for the Mellin transform we derive (cf. [2], Vol. II, Entry 2.16.3.15)

$$ S_{\mu,in} ( x) = {(2x)^{\mu+1} \over  \Gamma((1-\mu-in)/2) \Gamma((1-\mu+in)/2)}\int_0^\infty {t^{-\mu}\over t^2+ x^2}  K_{in}(t) dt.\eqno(2.35)$$
Hence with (2.34) we find

$$\left| S_{\mu,in} ( x)\right| \le   { A (2x)^{{\rm Re} \mu+1} \over   [\sinh(\pi n)]^{1/2} | \Gamma((1-\mu-in)/2) \Gamma((1-\mu+in)/2)|}\int_0^\infty {t^{-{\rm Re} \mu -1/4}\over t^2+ x^2}   dt$$

$$=  { A \ 2^{{\rm Re} \mu}  x^{-1/4} \Gamma((3/4- {\rm Re}\mu)/2) \Gamma((5/4+ {\rm Re}\mu)/2) \over   [\sinh(\pi n)]^{1/2} | \Gamma((1-\mu-in)/2) \Gamma((1-\mu+in)/2)|}.$$
\end{proof}

{\bf Theorem 5}. {\it   Let $-5/4 < {\rm Re}\mu < 0$ and a sequence $ a= \{a_n\}_{n\in \mathbb{N}} \in l_1$, i.e.   

$$||a||_{l_1}= \sum_{n=1}^\infty  |a_n|   < \infty.\eqno(2.36)$$
Then the discrete transformation $(1.5)$ can be inverted by the formula

$$a_n =  {2^{1-\mu} \over \pi^3} \sinh(\pi n)  \int_0^\infty   \Omega_n(x) f(x) dx,\ n \in \mathbb{N},\eqno(2.37)$$
where the kernel $\Omega_n(x)$ is defined by 

$$\Omega_n(x) = \int_0^\pi \sin\left(x\cosh(u)- {\pi\mu\over 2}\right)  \sinh(u) \sin(nu) du,\eqno(2.38)$$
and  integral  $(2.37)$ converges in the improper sense. }

\begin{proof}     Substituting the series (1.5) on the right-hand side of (2.37), we have

$$  {2^{1-\mu} \over \pi^3} \sinh(\pi n)  \int_0^\infty   \Omega_n(x) f(x) dx =  {2^{1-\mu} \over \pi^3} \sinh(\pi n) \lim_{T\to \infty}  \int_0^T   \Omega_n(x) $$

$$\times \sum_{m=1}^\infty a_m \  \Gamma\left({1-\mu-im\over 2}\right) \Gamma\left({1-\mu+im\over 2}\right) S_{\mu,im} (x) dx$$

$$= {2^{1-\mu} \over \pi^3} \sinh(\pi n) \lim_{T\to \infty}   \sum_{m=1}^\infty a_m \  \Gamma\left({1-\mu-im\over 2}\right) \Gamma\left({1-\mu+im\over 2}\right) $$

$$\times   \int_0^T   \Omega_n(x) S_{\mu,im} (x) dx,\eqno(2.39)$$
where the interchange of the order of integration and summation is due to condition (2.36) and inequality (2.33).  Then in order to pass to the limit under the series sign in (2.39), we estimate the remainder of the corresponding integral.  Indeed, recalling (2.35), (2.38) and the inequality for the modified Bessel function [4] $| K_{in}(x) |\le K_0(x)$, we deduce

$$\left| \Gamma\left({1-\mu-im\over 2}\right) \Gamma\left({1-\mu+im\over 2}\right) \int_T^\infty    \Omega_n(x) S_{\mu,im} (x) dx\right|  $$

$$=  \left|\int_T^\infty  (2x)^{\mu+1}   \int_0^\pi \sin\left(x\cosh(u)- {\pi\mu\over 2}\right) \right.$$

$$\times  \left.  \sinh(u) \sin(nu) du \int_0^\infty {t^{-\mu}\over t^2+ x^2}  K_{im}(t) dt dx\right| $$

$$\le 2^{{\rm Re} \mu+1} \pi \sinh(\pi) \left[ \left|\cos\left({\pi\mu\over 2}\right)\right|  + \left| \sin\left({\pi\mu\over 2}\right)\right| \right]$$

$$\times  \int_0^\infty  t^{-{\rm Re} \mu} K_{0}(t)  \int_T^\infty  {x^{{\rm Re} \mu+1}\over t^2+ x^2} dx dt $$

$$\le  2^{{\rm Re} \mu+1} \pi \sinh(\pi) \left[ \left|\cos\left({\pi\mu\over 2}\right)\right|  + \left| \sin\left({\pi\mu\over 2}\right)\right| \right]   {T^{{\rm Re} \mu} \over  - {\rm Re} \mu } \int_0^\infty  t^{-{\rm Re} \mu} K_{0}(t) dt, $$
and the latter expression tends to 0, when $T \to \infty$ by virtue of the condition ${\rm Re}\mu < 0$ and the convergence  of the latter integral (see (1.16), (1.17)). Therefore, passing to the limit under the series sign in (2.39), we appeal to (2.30) to obtain

$$   {2^{1-\mu} \over \pi^3} \sinh(\pi n)  \int_0^\infty   \Omega_n(x) f(x) dx = {2\over \pi} \sinh(\pi n)  \sum_{m=1}^\infty {a_m \over \sinh(\pi m) }$$

$$\times    \int_0^\pi  \sin(mu) \sin(nu) du = a_n.$$
Theorem 5 is proved.
\end{proof} 

A final result is the inversion theorem for the transformation (1.6). 

{\bf Theorem 6.} {\it Let ${\rm Re}\mu < 0$ and $f$ be a complex-valued function on $\mathbb{R}_+$ which is represented by the integral 

$$f(x) =   \int_{-\pi}^\pi \sin\left(x\cosh(u)- {\pi\mu\over 2}\right)   \varphi(u)  du,\quad x >0,\eqno(2.40)$$ 
where $ \varphi(u) = \psi(u)\sinh(u)$ and $\psi$ is a   $2\pi$-periodic function, satisfying the Lipschitz condition $(2.22)$ on $[-\pi, \pi]$.  Then for all $x >0$ the following inversion formula for  transformation $(1.6)$  holds

$$ f(x)  =   {2^{1-\mu}\over \pi^3} \sum_{n=1}^\infty  \sinh(\pi n)  \Omega_n (x) a_n,\eqno(2.41)$$
where $\Omega_n$ is defined by $(2.38)$.}

\begin{proof}    Substituting $f$ by formula (2.40) in (1.6), we invoke (2.30) and the definition of $\varphi$ to obtain after the interchange of the order of integration 

$$a_n = {2^{\mu} \pi^2 \over \sinh(\pi n)} \int_{-\pi}^\pi \psi (u) \sin(nu) du.\eqno(2.42)$$
This interchange is allowed due to the absolute convergence of the corresponding iterated integral via the asymptotic behavior (1.22) of the Lommel function at infinity for each fixed $n \in \mathbb{N}$.  Then,  substituting $a_n$ by (2.42) and $\Omega_n$ by (2.38) into the partial sum of the series (2.41) $S_N(x) $,   it becomes 

$$ S_N(x)  = {2\over \pi} \sum_{n=1}^N  \int_0^\pi \sin\left(x\cosh(t)- {\pi\mu\over 2}\right)  \sinh(t) \sin(nt) dt \int_{-\pi}^\pi \psi (u) \sin(nu) du.\eqno(2.43)$$
Hence,  calculating the sum via the known identity 
$$ \sum_{n=1}^N  \sin(nt) \sin(nu) = {1\over 4} \left[  {\sin \left((2N+1) (u-t)/2 \right)\over \sin( (u-t) /2)}  -  {\sin \left((2N+1) (u+t)/2 \right)\over \sin( (u+t) /2)} \right],$$
 we obtain from (2.43) 
 
 $$ S_N(x)  = {1\over 2\pi} \int_0^\pi \sin\left(x\cosh(t)- {\pi\mu\over 2}\right)  \sinh(t) dt \int_{-\pi}^\pi \left[ \psi (u) -\psi(-u)\right] $$
 
 $$\times   {\sin \left((2N+1) (u-t)/2 \right)\over \sin( (u-t) /2)} du.\eqno(2.44)$$
  Since $\psi$ is $2\pi$-periodic, we treat  the latter integral with respect to $u$ as follows 

$$  \int_{-\pi}^{\pi}  \left[ \psi (u) -\psi(-u)\right] \ {\sin \left((2N+1) (u-t)/2 \right)\over \sin( (u-t) /2)}  du $$

$$=  \int_{ t-\pi}^{t+ \pi} \left[ \psi (u) -\psi(-u)\right] \  {\sin \left((2N+1) (u-t)/2 \right)\over \sin( (u-t) /2)}  du $$

$$=  \int_{ -\pi}^{\pi}  \left[ \psi(u+t) - \psi(-u-t)\right] \  {\sin \left((2N+1) u/2 \right)\over \sin( u /2)}  du. $$
Moreover,

$$ {1\over 2\pi} \int_{ -\pi}^{\pi}  \left[ \psi(u+t)- \psi(-u-t) \right]  \  {\sin \left((2N+1) u/2 \right)\over \sin( u /2)}  du - \left[ \psi(t)- \psi(-t) \right] $$

$$=  {1\over 2\pi} \int_{ -\pi}^{\pi}  \left[ \psi(u+t)- \psi(t) + \psi (-t) - \psi(-u-t) \right]  \  {\sin \left((2N+1) u/2 \right)\over \sin( u /2)}  du.$$
When  $u+t > \pi$ or  $u+t < -\pi$ then we interpret  the value  $\psi(u+t)- \psi(t)$ by  formulas

$$\psi(u+t)- \psi(t) = \psi(u+t-2\pi)- \psi(t - 2\pi),$$ 

$$\psi(u+t)- \psi(t) = \psi(u+t+ 2\pi)- \psi(t +2\pi),$$ 
respectively.  Analogously, the value  $\psi(-u-t)- \psi(-t)$  can be treated.   Then   due to the Lipschitz condition (2.22) we have the uniform estimate
for any $t \in [-\pi,\pi]$

$${\left|  \psi(u+t)- \psi(t) + \psi (-t) - \psi(-u-t) \right| \over | \sin( u /2) |}  \le 2C \left| {u\over \sin( u /2)} \right|.$$
Therefore,  owing to the Riemann-Lebesgue lemma

$$\lim_{N\to \infty } {1\over 2\pi} \int_{ -\pi}^{\pi}  \left[ \psi(u+t)  - \psi(t) + \psi (-t) - \psi(-u-t)\right]  \  {\sin \left((2N+1) u/2 \right)\over \sin( u /2)}  du =  0\eqno(2.45)$$
for all $ t\in [-\pi,\pi].$    Besides, returning to (2.44), we estimate the iterated integral 
$$ \int_{0}^\pi \left|\sin\left(x\cosh(t)- {\pi\mu\over 2}\right) \right| \sinh(t)  \int_{ -\pi}^{\pi} \bigg\lvert \left[ \psi(u+t)  - \psi(t)+ \psi(-t) -\psi(-u-t)  \right]$$

$$\left.\times   {\sin \left((2N+1) u/2 \right)\over \sin( u /2)}  \right| du dt \le  2C [\cosh(\pi)-1]  $$

$$\times \left[\  \left|\cos\left({\pi\mu\over 2}\right)\right| +  \left|\sin\left({\pi\mu\over 2}\right) \right| \right]\int_{ -\pi}^{\pi}   \left| {u\over \sin( u /2)} \right| du < \infty.$$
 Consequently, via  the dominated convergence theorem it is possible to pass to the limit when $N \to \infty$ under the  integral sign, and recalling (2.45), we derive

$$  \lim_{N \to \infty}   {1\over 2\pi}  \int_{0}^\pi  \sin\left(x\cosh(t)- {\pi\mu\over 2}\right)  \sinh(t)   \int_{ -\pi}^{\pi}  \left[ \psi(u+t)  - \psi(t)+ \psi(-t) -\psi(-u-t)  \right] $$

$$\times  \  {\sin \left((2N+1) u/2 \right)\over \sin( u /2)}  du dt =  {1\over 2 \pi}  \int_{0}^\pi \sin\left(x\cosh(t)- {\pi\mu\over 2}\right)  \sinh(t)   $$

$$ \times \lim_{N \to \infty}  \int_{ -\pi}^{\pi}  \left[ \psi(u+t)-   \psi(t)  \right]  \  {\sin \left((2N+1) u/2 \right)\over \sin( u /2)}  du dt = 0.$$
Hence, combining with (2.44),  we obtain  by virtue of  the definition of $\varphi$ and $f$

$$ \lim_{N \to \infty}  S_N(x) =    \int_{0}^\pi  \sin\left(x\cosh(t)- {\pi\mu\over 2}\right)    \left[\varphi (t) +\varphi(-t) \right] dt = f(x),$$
where the integral (2.40) converges since $\varphi \in C[-\pi,\pi]$.  Thus we established  (2.41), completing the proof of Theorem 6.
 
\end{proof}

\bigskip
\centerline{{\bf Acknowledgments}}
\bigskip

\noindent The work was partially supported by CMUP, which is financed by national funds through FCT (Portugal)  under the project with reference UIDB/00144/2020.

\bigskip
\centerline{{\bf References}}
\bigskip
\baselineskip=12pt
\medskip
\begin{enumerate}

\item[{\bf 1.}\ ]  Yu.A. Brychkov, O.I. Marichev, N.V. Savischenko,   {\it Handbook of Mellin Transforms.}  Advances in Applied Mathematics, CRC Press,  Boca Raton, 2018.

\item[{\bf 2.}\ ] NIST Digital Library of Mathematical Functions. http://dlmf.nist.gov/, Release 1.0.17 of 2017-12-22. F. W. J. Olver, A. B. Olde Daalhuis, D. W. Lozier, B. I. Schneider, R. F. Boisvert, C. W. Clark, B. R. Miller and B. V. Saunders, eds.

\item[{\bf 3.}\ ] A.P. Prudnikov, Yu.A. Brychkov and O.I. Marichev, {\it Integrals and Series}. Vol. I: {\it Elementary
Functions}, Vol. II: {\it Special Functions}, Gordon and Breach, New York and London, 1986, Vol. III : {\it More special functions},  Gordon and Breach, New York and London,  1990.

\item[{\bf 4.}\ ] S. Yakubovich, {\it Index Transforms}, World Scientific Publishing Company, Singapore, New Jersey, London and
Hong Kong, 1996.

\end{enumerate}

\vspace{5mm}

\noindent S.Yakubovich\\
Department of  Mathematics,\\
Faculty of Sciences,\\
University of Porto,\\
Campo Alegre st., 687\\
4169-007 Porto\\
Portugal\\
E-Mail: syakubov@fc.up.pt\\

\end{document}